# Optimal Residential Demand Response Considering the Operational Constraints of Unbalanced Distribution Networks

Weiye Zheng, Wenchuan Wu, Boming Zhang
Department of Electrical Engineering
Tsinghua University, Beijing, China

Wanxing Sheng
China Electric Power Research Institute
Beijing, China

*Abstract*—As a typical approach of demand response (DR), direct load control (DLC) enables load service entity (LSE) to adjust electricity usage of home-end customers for peak shaving during DLC event. Households are connected in low voltage distribution networks, which is three phase unbalanced. However, existing works have not considered the network constraints and operational constraints of three phase unbalanced distribution systems, thus may ending up with decisions that deviate from reality or even infeasible in real world. This paper proposes residential DLC considering associated constraints of three phase unbalanced distribution networks. Numerical tests on a modified IEEE benchmark system demonstrate the effectiveness of the method.

## Nomenclature

The main symbols used throughout the paper are stated below for quick reference. Others are defined as required.

### A. Indices and Sets

| | |
|---|---|
| $i,j$ | Indices of buses, from 0 to $N$. |
| $\phi$ | Index of phases, typically phase $A$, $B$ and $C$. |
| $h$ | Index of households, from 1 to $H$. |
| $a$ | Index of appliances, from 1 to $A$. |
| $t$ | Index of time, from 1 to $T$. |
| $\mathrm{H}_i^\phi$ | Set of households at bus $i$, phase $\varphi$. |
| $\mathrm{A}_h$ | Set of appliances in household $h$. |
| $\mathrm{T}_{h,a}^{work}$ | Set of working time period of appliance $a$ in household $h$. |
| $\mathrm{T}_d$ | Set of DLC time period. |
| $N$ | Set of buses in distribution network. |
| $E$ | Set of lines in distribution network. |
| $G$ | Set of distributed generators. |

### B. Parameters

| | |
|---|---|
| $\eta_{h,a}$ | Power factor for appliance $a$ in household $h$. |
| $\mathbf{Z}_{ij} = \mathbf{R}_{ij} + j\mathbf{X}_{ij}$ | $3\times 3$ impedance matrix on line $ij$. |

### C. Variables

| | |
|---|---|
| $p_{h,a}(t), q_{h,a}(t),$ $s_{h,a}(t)$ | Active, reactive and complex power of appliance $a$ in household $h$ in time $t$ respectively. |
| $p_{Gi}^\phi, q_{Gi}^\phi$ | Active and reactive power output of distributed generator at bus $i$, phase $\varphi$, respectively. |
| $T_h^{in}(t)$ | Indoor temperature of household $h$ in time $t$. |
| $\mathbf{v}_i(t)$ | $3\times 1$ vector of squared voltage magnitude for three phases at bus $i$ at time $t$. |
| $\mathbf{P}_{ij}(t), \mathbf{Q}_{ij}(t)$ | $3\times 1$ vector of active and reactive branch flow for three phases on line $ij$ at time $t$ respectively. |
| $\mathbf{p}_i(t), \mathbf{q}_i(t)$ | $3\times 1$ vector of net injected active and reactive power for three phases at bus $i$ at time $t$ respectively. |

## I. Introduction

DEMAND side management is a mechanism to help reduce peak load and adjust elastic load to improve system efficiency [1]. Most of previous works consider the supply-demand balance in demand response (DR) in an abstract way and neglect network constraints [1-5]. However, households are connected in a distribution network that follows power flow constraints and operational constraints. Therefore, the optimal decisions provided by these schemes may violate some operational constraints and thus are infeasible in real system. Distribution network constraints are first incorporated in DR in [6], but it has the underlying assumption that both the distribution network and the households are three phase balanced. In reality, households are connected by a low voltage (LV) distribution network in most cases, which is typically three phase unbalanced [7]. Also, one of the most common forms is single phase customer [8].

None of previous works have considered the underlying three-phase unbalanced distribution networks when performing DR. This paper extends the residential demand response model to a more elaborate and realistic formulation by considering three phase unbalanced distribution network. We consider a residential direct load control (DLC) program where customers give the load service entity (LSE) the option to remotely shut down appliances for a certain period of time (e.g., peak hours) by signing up a contract in advance [9, 10]. In a DLC event, the home energy management systems (HEMS) follow the control signals from the LSE and schedule their appliance consumption to fulfill the following DLC objectives: (i) social welfare maximization (i.e., the customer utilities minus LSE cost); (ii) system demand under operational limit during peak hours; (ii) the appliance operational constraints, power flow constraints, network security and system operational constraints are satisfied.

The remainder of the paper is organized as follows. In Section II describes problem formulation. Section III details the results of several numerical tests to investigate the effectiveness of the proposed DLC scheme. Section IV concludes the paper.

## II. DLC Problem Formulation

Residential DLC over a three phase unbalanced distribution network operated by an LSE is considered in this section. Each load bus is connected with a set of households and there are $H$ households $H = \{h_1, h_2, ..., h_H\}$ in the system. In each household, an HEMS is managing the appliances $A_h = \{a_{h,1}, a_{h,2}, ..., a_{h,A}\}$ such as air conditioners (ACs), dryers, washers, lighting, etc. Suppose there is a two-way communication between the LSE and the household [11]. In a contract-based DLC program, the HEMS receive DLC event from LSE and shed or shift their demands to fulfill the objective.

### A. Appliance Model

Conventionally appliances in DLC can be classified into critical, interruptible and deferrable loads [12]. Thermostatically controlled appliances are usually a kind of interruptible loads, but they are specified into a single category in this paper due to their unique features.

Customer utility function of different types of appliances and their operational constraints are described as follow:

- For interruptible loads such as optional lighting and plug loads, the utility depends on the power it consumes at time $t$ and can be written in the form of:

$$U_{h,a}(p_{h,a}) = \sum_{t=1}^{T} U_{h,a}(p_{h,a}(t)). \qquad (1)$$

- For deferrable loads such as washers and dryers, the utility depends on the total energy consumption and may be time sensitive:

$$U_{h,a}(p_{h,a}) = U_{h,a}(\sum_{t=1}^{T} p_{h,a}(t)\Delta t, t). \qquad (2)$$

Also, the cumulative energy consumption must lie within an acceptable interval:

$$\underline{E}_{h,a} \leq \sum_{t \in T} p_{h,a}(t)\Delta t \leq \overline{E}_{h,a}. \qquad (3)$$

- For thermostatically controlled appliances such as ACs and heaters, the customer would be happier if the indoor temperature is closer to the most comfort temperature $T_h^{comf}(t)$. By denoting a demand vector $p_{h,a} = \{p_{h,a}(t), t = 1, 2, ..., T\}$, the utility can be defined in a general form of:

$$U_{h,a}(p_{h,a}) = \sum_{t=1}^{T} U_{h,a}(T_h^{in}(t), T_h^{comf}(t)) \qquad (4)$$

Also, the temperature evolution equation must hold [6]:

$$T_h^{in}(t) = T_h^{in}(t-1) + \alpha_h(T_h^{out}(t) - T_h^{in}(t-1)) + \beta_h p_{h,a}(t) \qquad (5)$$

The indoor temperature should be kept within a bearable range:

$$\underline{T}^{in}(t) \leq T^{in}(t) \leq \overline{T}^{in}(t) \qquad (6)$$

- Besides, the following common constraints exist in all above appliances.

The power factor constraint:

$$\eta_{h,a} = \frac{p_{h,a}(t)}{|s_{h,a}(t)|} \qquad (7)$$

Operational limit of appliances:

$$\underline{p}_{h,a}(t) \leq p_{h,a}(t) \leq \overline{p}_{h,a}(t) \qquad (8)$$

Non-working hour constraints:

$$\underline{p}_{h,a}(t) = \overline{p}_{h,a}(t) = 0, \qquad \forall t \notin T_{h,a}^{work} \qquad (9)$$

- Critical loads such as refrigerators, cooking and critical lighting are modeled as given loads that have to be satisfied and should not be shifted or shed at any time.

### B. Distributed Generation (DG) Model

Here we consider the integration of photovoltaic (PV) power generators into distribution network. Suppose that the output of a DG follows the signal from maximum power point tracking (MPPT). Therefore, all DGs generate the maximum allowable active power, and the reactive power can be controllable and optimized. We have

$$p_{Gi}(t) = \overline{p}_{Gi}(t). \qquad (10)$$

$$\underline{q}_{Gi}(t) \leq q_{Gi}(t) \leq \overline{q}_{Gi}(t) \qquad (11)$$

### C. Distribution Network Model

A distribution network is modeled as a connected graph (N, E). Bus 0 denotes the root or the point of common coupling (PCC). Specifically, the feeder voltage $V_0$ is fixed at reference voltage regulated by substation side [13]. And $p_0(t)$ is the power injected to the distribution system. Each load bus $i \in N \setminus \{0\}$ provides power to a set of households connected to the bus. The net injected power at bus $i$ phase $\phi$ satisfies:

$$p_i^\phi(t) = p_{i,crit}^\phi(t) + \sum_{h \in H_i^\phi, a \in A_h} p_{h,a}(t) - p_{Gi}^\phi(t) \qquad (12)$$

$$q_i^\phi(t) = q_{i,crit}^\phi(t) + \sum_{h \in H_i^\phi, a \in A_h} q_{h,a}(t) - q_{Gi}^\phi(t) \qquad (13)$$

where $p_{i,crit}^\phi, q_{i,crit}^\phi$ denote critical loads.

For distribution system, Distflow branch equations are typically used to describe power flow. Extending the original formulation in [14] to three phase unbalanced cases, we have:

$$\sum_{i:i \to j}[P_{ij} - P_{ij}^{loss}(P,Q)] + p_j = \sum_{k:j \to k} P_{jk}, \quad j \in N$$

$$\sum_{i:i \to j}[Q_{ij} - Q_{ij}^{loss}(P,Q)] + q_j = \sum_{k:j \to k} Q_{jk}, \quad j \in N \qquad (14)$$

$$v_j = v_i - 2(R_{ij}P_{ij} + X_{ij}Q_{ij}) + \Delta v_{ij}(P,Q) \qquad (15)$$

$$v_0^\phi = V_{ref}^2, \quad \forall \phi \qquad (16)$$

where $P_{ij} \coloneqq \{P_{ij}^\phi(t), \forall \phi, t\}$, $Q_{ij} \coloneqq \{Q_{ij}^\phi(t), \forall \phi, t\}$, $p_j \coloneqq \{p_j^\phi(t), \forall \phi, t\}$, $q_j \coloneqq \{q_j^\phi(t), \forall \phi, t\}$ and $v_j \coloneqq \{v_j^\phi(t), \forall \phi, t\}$. Although $P_{ij}^{loss}(P,Q), Q_{ij}^{loss}(P,Q), \Delta v_{ij}(P,Q)$ are all non-linear terms, we can extend the Taylor expansion technique in our previous work [15] to three phase unbalanced cases and ensure (14) and (15) to be linear constraints. Detailed derivation is available in [16]. Elaborate expression is attached in Appendix.

Since distribution network has been formulated in details, now the security constraints of branch flow and voltage magnitude can be considered in the model now:

$$\underline{P}_{ij}^{\phi} \leq P_{ij}^{\phi}(t) \leq \overline{P}_{ij}^{\phi}$$
$$\underline{Q}_{ij}^{\phi} \leq Q_{ij}^{\phi}(t) \leq \overline{Q}_{ij}^{\phi} \quad (17)$$
$$\underline{V}_i^2 \leq v_i^{\phi}(t) \leq \overline{V}_i^2 \quad (18)$$

### D. DLC Model

Given the DLC event, the system demand constraint can be formulated as:

$$|s_0(t)| = \sqrt{\|\mathbf{p}_{G0}(t)\|_2^2 + \|\mathbf{q}_{G0}(t)\|_2^2} \leq \overline{S}, \quad \forall t \in T_d \quad (19)$$

where $\mathbf{p}_{G0}(t), \mathbf{q}_{G0}(t)$ is $3 \times 1$ total active and reactive power injection vector at PCC and it is given by:

$$\mathbf{p}_{G0}(t) = \sum_{j:0 \to j} \mathbf{P}_{0j}(t) \quad (20)$$

$$\mathbf{q}_{G0}(t) = \sum_{j:0 \to j} \mathbf{Q}_{0j}(t) \quad (21)$$

The objective is to find a set of optimal demand and generation vectors to maximize the aggregate utilities of the households and minimize the overall cost of LSE in distribution network:

$$\max_{p,q,v} \quad \kappa \sum_{h \in H} \sum_{a \in A_h} U_{h,a}(\mathbf{p}_{h,a}) - \sum_{t \in T} \sum_{j \in G \cup \{0\}} C_t(p_{Gj}(t)) \quad (22)$$

$$s.t. \quad (3), (6)-(20)$$

## III. SIMULATION RESULTS

The optimization package CPLEX [17] was used to solve the model in conjunction with Matlab R2015a, using a machine with an Intel Core i5-3210M 2.50-GHz processor and 8 GB of RAM. The proposed scheme is tested on a 10-bus distribution system derived from IEEE 33-bus distribution system [18], as shown in Fig. 1. The cost function of LSE is in quadratic form:

$$C_t(p_{Gi}(t)) = a(t)[p_{Gi}(t)]^2 + b(t)p_{Gi}(t) + c(t) \quad (23)$$

which includes both generation cost of PVs and the cost of power acquired at PCC from external networks.

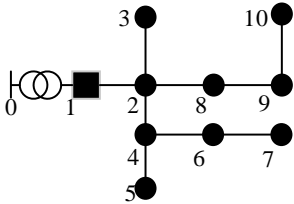

Fig. 1    The modified IEEE standard distribution system

15 households are connected to each of load bus 2 ~10. Time interval is 1 hour. 4 PVs are located at bus {3, 6, 8, 9}. All households and PVs are randomly distributed at phase A, B and C. 5 types of appliances including ACs, washers, dryers, optional lighting and plug loads are taken into account in simulation. $\eta_{h,a}$ is randomly selected from [0.8, 0.9].

(1) ACs (thermostatically controlled load): In simulation, $\alpha = 0.9$, and $\beta$ is randomly chosen from [-0.008, -0.005]. Bearable temperature range is set as [70F, 79F], and the most comfort temperature is randomly chosen from [73F, 76F]. The upper bound and lower bound of power are $\overline{p}_{h,a} = 3.5kW, \underline{p}_{h,a} = 0kW$, respectively. The outside temperature of a day is depicted in Fig. 2, which is a typical summer day in Southern California [6]. The utility takes the form of

$$U_{h,a}(p_{h,a}) = -\sum_{t=1}^{T} b_{h,a}(T_h^{in}(t) - T_h^{comf}(t))^2 \quad (24)$$

where $b_{h,a}$ is positive constant.

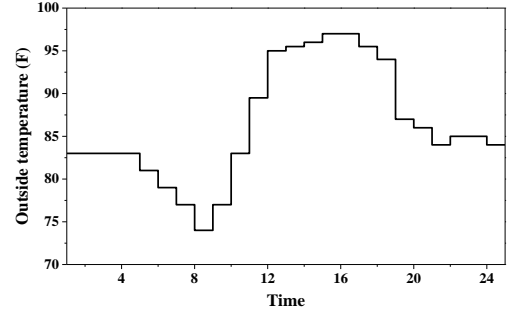

Fig. 2    Outside temperature of a day

(2) Washers (deferrable load): The time when customers arrive at their home $t_{h,arrive}$ is randomly selected from [17, 19]. Starting time of the washer is an hour after $t_{h,arrive}$ and it must finish its task within 2 hours. The bound of power is set as $\overline{p}_{h,a} = 700W, \underline{p}_{h,a} = 0W$. Total energy consumption is set as $\overline{E}_{h,a} = \underline{E}_{h,a} = 900Wh$. The utility takes the form of:

$$U_{h,a}(p_{h,a}) = b_{h,a}(\sum_{t=1}^{T} p_{h,a}(t)\Delta t)$$
$$-\sum_{t=1}^{T} t[p_{h,a}(t) - p_{h,a}^{pref}(t)]^2 \quad . \quad (25)$$

(3) Dryers (deferrable load): The dryer starts working at $t_{h,arrive} + 3$ and must finish before $t = 1$. The bound of active power is set as $\overline{p}_{h,a} = 5kW, \underline{p}_{h,a} = 0W$, and cumulative energy consumption range is set as $\overline{E}_{h,a} = 9kWh, \underline{E}_{h,a} = 4.5kWh$. The utility is in the same form as that of washers.

(4) Optional lighting (interruptible load): It works from 7 p.m. to 7 a.m. The maximum and minimum powers are $\overline{p}_{h,a} = 1kW, \underline{p}_{h,a} = 0.5kW$. The utility function is in the form of

$$U_{h,a}(p_{h,a}) = -b_{h,a}(p_{h,a}(t) - p_{h,a}^{pref}(t))^2 \quad . \quad (26)$$

(5) Plug loads (interruptible load): The maximum and minimum powers are $\overline{p}_{h,a} = 500W, \underline{p}_{h,a} = 0W$. The utility function takes the same form as that of optional lighting.

The reference voltage is set as $V_{ref} = 4.16kV$, while the minimum allowed voltage $\underline{V}_i = 4.05kV$ [19]. Normalized curves of critical load, PV active power output and cost of LSE are shown in Fig. 3 [15, 20].

In simulation, the proposed DLC scheme is compared with the baseline without DLC (W/O DLC in short) and conventional DLC.

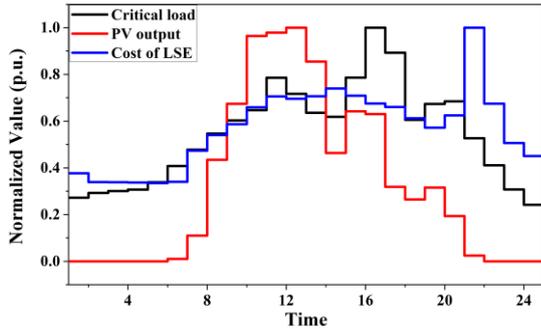

Fig. 3　　Normalized curves of critical load, PV active power output and generation price.

**W/O DLC** is defined as follows:
- Only consider customer utilities maximization
- The AC keeps the indoor temperature to the most comfortable temperature all day.
- The washer and dryer run at their maximum power until the maximum energy requirement is met.
- Optional lighting and plug loads use the power as they request.
- The curtailment and security constraints in direct load control model are not considered.

**Conventional DLC** considers supply-demand in an abstract way. The distribution network model in Section II. C is replaced by a conventional overall supply-demand balance constraint. Please note that voltage security constraint cannot be considered for there is no network model. Other remains the same as the proposed scheme.

The optimal schedule of DGs and households generated by all three approaches is simulated over real distribution system using three phase modelling to get the actual power flow.

For the base case without DLC, the load profile at the root $s_0(t)$ is drawn by dotted line in black in Fig. 4. There is a peak at $t = 22$. The dotted lines in Fig. 5 shows the minimum bus voltage magnitude in the distribution network over time for all three phases for base case. It can be observed that bus voltage violates the security constraints during the peak. The voltage drop is more severe when the demand goes higher. This relationship can be explained by the voltage equation (15).

To simulate a DLC event, the demand limit $\bar{s}$ is set as 0.95MVA during time period [19, 24]. The simulated load profiles of PCC $s_0(t)$ with conventional DLC and the proposed DLC are depicted by the solid line in Fig. 4 in red and blue, respectively. Both methods can effectively mitigate the peak and successfully keep the system demand under the limit during the DLC event. But please note that the demand limit of conventional DLC is not strictly satisfied and it was 0.955MVA at $t = 22$. The reason is that conventional DLC neglects power losses.

Although conventional DLC can roughly manage demand, its decision violates voltage security constraints and, as a result, the decision is infeasible in real system. As shown by the blue solid line in Fig. 5, the voltage violation in phase C at $t = 22$ cannot be relieved by conventional DLC.

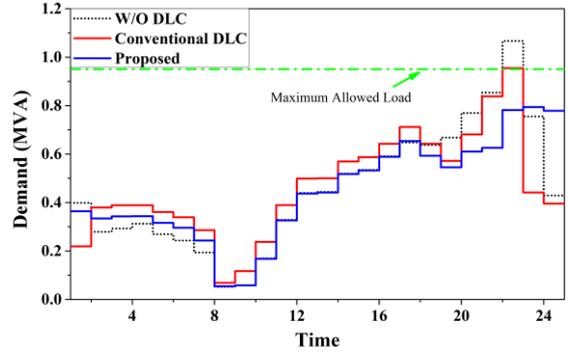

Fig. 4　　Comparison of load profile of PCC $s_0(t)$ without DLC, with conventional DLC and with the proposed scheme.

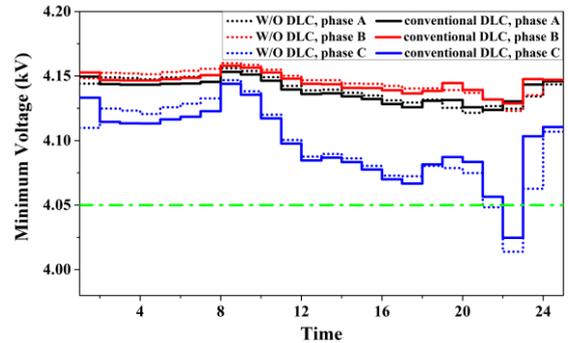

Fig. 5　　Minimum bus voltage profiles without and with conventional DLC.

On the contrary, by considering three phase unbalanced distribution network constraint, the proposed DLC can not only shave the peak, but also provide a feasible decision that satisfies physical and security constraints in real system. Detailed voltage profiles for all three phases of the proposed scheme are depicted in solid lines in Fig. 6 in black, red and blue, respectively. The voltage magnitudes lie within operational limit.

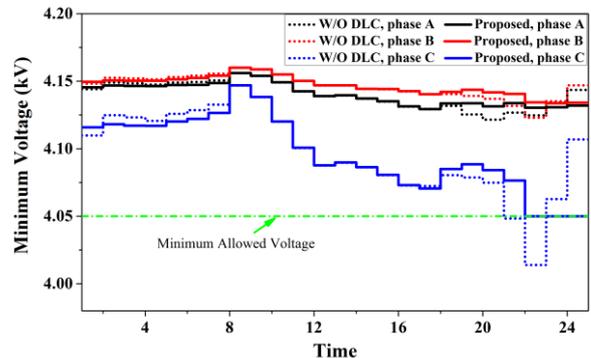

Fig. 6　　Minimum bus voltage profiles without and with the proposed DLC.

To investigate further how the proposed scheme works, Fig. 7 shows the load profile of the appliances in one of the households at bus 8 in phase C without DLC and with the proposed DLC. Both the washer and the dryer have been shifted and shaved, and the total energy consumption of the dryer is also reduced. Interruptible loads have also been curtailed to some extent. This can, in part, explain the peak shaving in Fig. 4. With the proposed DLC scheme, the cost is reduced from 160.4 to 145.3 for the sake of peak shaving. Due to limit of space, test results on a larger system are included in supplementary materials, which is available online [21].

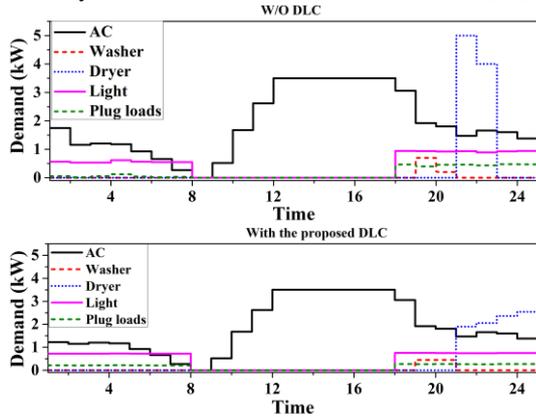

Fig. 7  Load profile of the appliances in one of the households at bus 8 in phase C without DLC and with the proposed DLC.

IV. CONCLUSION

We described residential direct load control considering the underlying three phase unbalanced distribution network constraints, which is joint maximization of LSE and households with consideration of their operational constraints. Test results on a modified IEEE distribution system show that the proposed scheme can effectively shave peak in unbalanced distribution network. Although conventional DLC considering abstract supply-demand matching can also manage demand in a rough fashion, its decision is not feasible in real system because it violates the physical and security constraints. Since both cost structure of LSE and customers' utility functions are private information, future works includes designing efficient decentralized scheme to address the privacy issues.

APPENDIX

Elaborate form of the nonlinear terms in (14)-(15):

$$\boldsymbol{P}_{ij}^{loss}(P,Q) = \boldsymbol{P}_{ij}.*(\hat{\boldsymbol{r}}_{ik}\boldsymbol{P}_{ik} + \hat{\boldsymbol{x}}_{ik}\boldsymbol{Q}_{ik}) + \boldsymbol{Q}_{ik}.*(\hat{\boldsymbol{r}}_{ik}\boldsymbol{Q}_{ik} - \hat{\boldsymbol{x}}_{ik}\boldsymbol{P}_{ik}) \quad (23)$$

$$\boldsymbol{Q}_{ij}^{loss}(P,Q) = \boldsymbol{P}_{ij}.*(\hat{\boldsymbol{x}}_{ik}\boldsymbol{P}_{ik} - \hat{\boldsymbol{r}}_{ik}\boldsymbol{Q}_{ik}) + \boldsymbol{Q}_{ik}.*(\hat{\boldsymbol{r}}_{ik}\boldsymbol{P}_{ik} + \hat{\boldsymbol{x}}_{ik}\boldsymbol{Q}_{ik}) \quad (24)$$

$$\Delta \boldsymbol{v}_{ij}(P,Q) = [\boldsymbol{Z}_{ik}(\boldsymbol{S}_{ik}^{*}./\boldsymbol{v}_{i}^{*})].*[\boldsymbol{Z}_{ik}^{*}(\boldsymbol{S}_{ik}./\boldsymbol{v}_{i})] \quad (25)$$

where ".*" and "./" denote element-wise multiplication and division, respectively.

Due to limit of space, detailed expression of Taylor expansion is not included in this paper. For readers who are interested, please refer to the Appendix in [16].